\documentclass[12pt]{amsart}
\usepackage{amsmath}\usepackage{amssymb}
\usepackage{amscd}
\usepackage[all]{xy}

\def\lra{\longrightarrow}
\def\map#1{\,{\buildrel #1 \over \lra}\,}

\def\duBois{du\,Bois\ }
\def\zar{{\mathrm{zar}}}
\def\cdh{{\mathrm{cdh}}}

\def\cO{\mathcal O}
\newcommand{\mfc}{\mathfrak{c}}
\newcommand{\mfm}{\mathfrak{m}}
\newcommand{\bbH}{\mathbb H}

\newcommand{\bbHcdh}{\bbH_{\cdh}}

\newcommand{\R}{\mathbb{R}}
\newcommand{\Q}{\mathbb{Q}}
\newcommand{\Z}{\mathbb{Z}}

\def\Spec{\operatorname{Spec}}

\numberwithin{equation}{section}

\theoremstyle{plain}
\newtheorem{thm}[equation]{Theorem}
\newtheorem{prop}[equation]{Proposition}
\newtheorem{lem}[equation]{Lemma}
\newtheorem{cor}[equation]{Corollary}

\theoremstyle{definition}
\newtheorem{defn}[equation]{Definition}
\newtheorem{ex}[equation]{Example}

\newtheorem{Wk-mod}[equation]{$W_\Gamma(k)$-module structures}

\theoremstyle{remark}
\newtheorem{rem}[equation]{Remark}

\newtheorem{question}[equation]{Question}
\begin {document}

\title[]
{$K_2$-regularity and normality}  
\date{\today}
\author{Christian Haesemeyer}
\address{School of Mathematics and Statistics, University of Melbourne,
VIC 3010, Australia}
\email{christian.haesemeyer@unimelb.edu.au}\urladdr{https://blogs.unimelb.edu.au/christian-haesemeyer}
\thanks{Haesemeyer was supported by ARC DP-210103397}

\author{Charles A.\,Weibel}
\address{Math.\ Dept., Rutgers University, New Brunswick, NJ 08901, USA}
\email{weibel@math.rutgers.edu}\urladdr{http://math.rutgers.edu/~weibel}
\thanks{Weibel was supported by NSF grant 2001417}


\begin{abstract}
We take a fresh look at the relationship between $K$-regularity and
regularity of schemes, proving two results in this direction. First, we
show that $K_2$-regular affine algebras over fields of characteristic
zero are normal.  Second, we improve on Vorst's $K$-regularity bound
in the case of local complete intersections; this is related to recent
work on higher du Bois singularities.
\end{abstract}
\maketitle

Recall that a ring $R$ is {\it $K_n$-regular} if and only if $NK_s(R) = 0$
for all $s\leq n$,
where $NK_n(R)$ is $K_n(R[t])/K_n(R)$.
In his 1979 paper \cite{Vorst1}, Ton Vorst conjectured that an affine
$k$-algebra $A$ of dimension $d$ is regular if, and only if, it is
$K_{d+1}$-regular. 
This conjecture was proved in characteristic zero
by the authors and Corti\~nas in \cite{chw-vorst}, and by Kerz,
Strunk, and Tamme over perfect fields of characteristic $p$ in
\cite{KST21}, drawing on previous work of Geisser and Hesselholt in
\cite{GH12}. 

In this paper, we make two improvements on existing results. 
In the first (Theorem \ref{thm:K_2-regular-normal}), we generalize 
Vorst's original result in a different direction by proving 
that $K_2$-regularity implies normality in arbitrary dimensions: 

\begin{thm}\label{thm A}
Let $A$ be a noetherian commutative ring containing $\Q$ with
finite normalization. If $A$ is $K_2$-regular then $A$ is normal.
\end{thm}

\begin{rem}
	The converse of Theorem \ref{thm A} is false in
	dimensions larger than $1$. Indeed, \cite{W2001} provides numerous
	examples of normal surface singularities that are not even
	$K_{-1}$-regular, let alone $K_2$-regular.
\end{rem}

In the second (Theorem  \ref{thm:Vorst-lci}), we  
begin investigating
whether Vorst's bound is the best possible. Combining the $K$-theoretic
calculations from \cite{chww-bass} and algebraic-geometric results
from \cite{MP22}, we show 
that much sharper bounds are possible,
at least in the case of local complete intersections: 

\begin{thm}\label{thm B}
Let $X\in\mathrm{Sch}/k$ be an affine local complete intersection 
scheme over a field of characteristic zero, and $p\geq 1$. If $X$ is
$K_{p+1}$-regular, then $X$ is regular in codimension $2p$. In particular, if $X$ has dimension $d$ and is $K_q$-regular for some 
$q \geq d/2 + 1$, then $X$ is regular.
\end{thm}

\begin{ex}
Consider the case of a surface $X$. Theorem \ref{thm A} says that if
$X$ is $K_2$-regular, then it is normal. Theorem \ref{thm B} says that
if $X$ is moreover a local complete intersection, then it is, in fact,
non-singular. In contrast, the main result of \cite{chw-vorst}
requires $K_3$-regularity for the latter conclusion.
\end{ex}

\begin{rem}
The bound in Theorem \ref{thm B} can be sharpened, see Corollary \ref{cor:dB-codim} and Remark \ref{rem:Vorst-lci}. 	
\end{rem}


\section{$K_2$-regularity implies normality.}\label{sec:normality}
In this section we prove that a $K_2$-regular algebra of finite type
over a field of characteristic zero is necessarily normal.
This generalizes the arguments that Vorst used in \cite{Vorst1} 
to prove that $K_2$-regular curves over fields must be 
regular.
By \cite[1.9]{Vorst1}, a commutative ring $A$ is $K_n$-regular
iff $A_{\mfm}$ is $K_n$-regular for every 
maximal ideal $\mfm$ of $A$.
This reduces the question to the case in which $A$ is local.

Throughout this section, $A$ will be a reduced noetherian local ring,
with nonzero maximal ideal $\mfm$, 
and $B$ is a finite birational extension of $A$.  An important case
is when the normalization of $A$ is finite over $A$. 
This is the case for example when
$A$ is a localization of some affine algebra
over a field. We will write $\mfc$ 
for the conductor of $A\to B$, that is, 
the annihilator of the finitely generated $A$-module $B/A$.

It is well-known that there is a Milnor square of rings
\begin{equation}\label{square:conductor}
\xymatrix{ A\ar[d] \ar[r]& B \ar[d] \\ A/\mfc \ar[r] & B/\mfc.}
\end{equation}
From this we obtain a diagram of relative $K$-theory exact sequences,
for $n\in\Z$,
\begin{equation}\label{ladder}
\xymatrix{
NK_{n+1}(A/\mfc) \ar[d] \ar[r]&
NK_n(A,\mfc) \ar[d] \ar[r]& NK_n(A) \ar[d] \ar[r]& NK_n(A/\mfc) \ar[d] 
\\
NK_{n+1}(B/\mfc) \ar[r]&
NK_n(B,\mfc) \ar[r]& NK_n(B) \ar[r]& NK_n(B/\mfc). 
	}\end{equation}
	
If $R$ is a commutative ring, we shall write $NU(R)$  
	for the group of units $R[t]^\times/R^\times$; it is 
	in 1--1 correspondence with the set $\mathrm{Nil}(R)[t]$.
	It is also the top of the $\gamma$-filtration of $NK_1(R)$.
	
\begin{lem}\label{lem:K_1-regular_preserved}
Suppose the conductor $\mfc$ is $\mfm$-primary.\  
If $A$ is $K_1$-regular, then so is $B$.
\end{lem}

\begin{proof}
By hypothesis, $A/\mfc$ and $B/\mfc$ are Artinian;
therefore $NK_n(A/\mfc)\!= NK_n(B/\mfc) = 0$ for $n\leq0$. 
By excision, the maps $NK_n(A,\mfc)\to
NK_n(B,\mfc)$ are isomorphisms for $n\leq 0$ and onto for $n =1$. 
It now follows from (\ref{ladder}) and the $5$-Lemma that
$NK_n(A)\to NK_n(B)$ is an isomorphism for $n <0$ and onto for $n=0$;
therefore, $B$ is $K_0$-regular.

Since $NK_1(A) = 0$, the map $NK_2(A/\mfc)\to
NK_1(A,\mfc)$ is onto. Because $NK_1(A,\mfc)\to
NK_1(B,\mfc)$ is onto by excision, it follows that
$NK_2(B/\mfc)\to NK_1(B,\mfc)$ is onto as well. Using
exactness of the bottom row of (\ref{ladder}), we conclude that
$NK_1(B)\to NK_1(B/\mfc)$ is injective. 

Since $B/\mfc$ is Artinian, $NK_1(B/\mfc) = NU(B/\mfc).$
Since $NK_1(B)$ is a subgroup of $NK_1(B/\mfc)$,
it follows that $NK_1(B) = NU(B)$ as well. 
But $B$ is birational over $A$, and in particular
is reduced, so $NU(B) = 0$. This shows that 
$NK_n(B) = 0$ for all $n\leq 1$, that is, $B$ is $K_1$-regular.	
\end{proof}

\begin{cor}\label{cor:K_1-regular_radical_conductor}
Suppose the conductor $\mfc$ is $\mfm$-primary.  
If $A$ is $K_1$-regular then $\mfc = \mfm$ and
$B/\mfc$ is a product of field extensions of $A/\mfm$.
\end{cor}

\begin{proof}
By Lemma \ref{lem:K_1-regular_preserved}, $NK_1(A) = NK_1(B) = NK_0(A)
= NK_0(B) = 0$. Because $NK_0(A,\mfc)\to NK_0(B,\mfc)$
is an isomorphism by excision, (\ref{ladder}) implies that
$NK_1(A/\mfc)\to NK_1(B/\mfc)$ is also an
isomorphism. In particular, we have an isomorphism
$NU(A/\mfc)\map{\cong} NU(B/\mfc)$ and hence an isomorphism
$\mathrm{Nil}(A/\mfc)\map{\cong} \mathrm{Nil}(B/\mfc)$ 
of nil radicals. 

Let $x\in B$ be an element such that the reduction
$\overline{x}$ of $x$ modulo $\mfc$ is in
$\mathrm{Nil}(A/\mfc) \cong \mathrm{Nil}(B/\mfc)$. Then
$\overline{x}\left(B/\mfc\right)\subseteq A/\mfc$, and
therefore $xB\subseteq A + \mfc = A$. Thus, $x\in\mfc$
by definition of the conductor, so that $\overline{x} = 0$. 
It follows that $\mathrm{Nil}(A/\mfc) = \mathrm{Nil}(B/\mfc) =0$. 
Because $A/\mfc$ and $B/\mfc$ are Artinian rings,
this implies that $\mfc = \mathfrak{m}$, and that $B/\mathfrak{m} = \prod_i L_i$ is a finite product of field extensions of $K = A/\mathfrak{m}$.
\end{proof}

We will write $K$ for $A/\mfm$ and $L$ for 
$B/\mathfrak{m}= \prod_i L_i$.

\begin{lem}\label{KABI}
$NK_1(A,B,\mfm)$ is zero iff $L/K$ is separable.
\end{lem}

\begin{proof}
Now $NK_1(A,B,\mfm)$ is a direct summand of
$K_1(A[t],B[t],\mfm[t])$,
and the latter group is computed in \cite[Theorem 0.2]{GW83}
as
\[
K_1(A[t],B[t],\mfm[t]) \cong~ 
\mfm/\mfm^2\otimes_{L} \Omega_{L[t]/K[t]} 
\cong~ \mfm/\mfm^2\otimes_{L} \Omega_{L/K}\otimes_L L[t]
\]
It is standard that $L/K$ is separable if and only if $\Omega_{L/K}=0$,
so $K_1(A[t],B[t],\mfm[t])=0$ if and only if $L/K$ is separable.
\end{proof}

\begin{cor}\label{L/Kseparable}
If $A$ is $K_1$-regular and
$L/K$ is not separable, then $NK_2(B)\ne0$.
\end{cor}

\begin{proof}
We refer to diagram \eqref{ladder}. Since $A$ is $K_1$-regular,
$NK_1(A,\mfm)=0$. Because $L/K$ is a product of fields, 
$NK_2(B)\cong NK_2(B,\mfm)$.
From the exact sequence $NK_2(B,\mfm)\to NK_1(A,B,\mfm)\to NK_1(A,\mfm)$,
we see that $NK_2(B)$ maps onto $NK_1(A,B,\mfm)$, 
which is nonzero by Lemma \ref{KABI}.
\end{proof}

\begin{cor}\label{cor:K_2-regular_preserved}
Suppose that the conductor $\mfc$ is $\mathfrak{m}$-primary and $A$ is $K_2$-regular. 
Then $B$ is $K_2$-regular iff $B/\mfc$ is a product of finite separable field extensions of $A/\mfm = A/\mfc$.
\end{cor}

\begin{proof}
Suppose first that $B$ is $K_2$-regular. By Corollary \ref{cor:K_1-regular_radical_conductor}, $\mfc = \mfm$ 
and $B/\mfc$ is a product of finite field extensions of $A/\mfm$. Corollary \ref{L/Kseparable} implies that the extensions are separable.  

Now suppose that $B/\mfc$ is a product of finite separable field extensions of $A/\mfc$; by Lemma \ref{lem:K_1-regular_preserved}, $B$ is $K_1$-regular. 
It remains to prove that $NK_2(B) = 0$. By Corollary
\ref{cor:K_1-regular_radical_conductor}, $A/\mfc$ and
$B/\mfc$ are regular rings, so $NK_n(A/\mfc) =
NK_n(B/\mfc) = 0$ for all $n\in\Z$. Using
(\ref{ladder}), it follows that $NK_n(A) \cong
NK_n(A,\mfc)$ and $NK_n(B) \cong NK_n(B,\mfc)$ for all integers $n$.
Lemma \ref{lem:K_1-regular_preserved} implies that
$NK_1(A,\mfc) = NK_1(B,\mfc) = 0$.
It follows that the following sequence is exact:
\[NK_2(A,\mfc) \to NK_2(B,\mfc) \to NK_1(A,B,\mfc) \to 0.\]
Now Lemma \ref{KABI} implies the assertion. 
\end{proof}

We now treat a special case of the main result of this section. 

\begin{prop}\label{prop:isolated-normal}
Suppose that $\left(A, \mfm\right)$ is a noetherian local ring,
$B/A$ is a finite birational extension, and that the conductor
$\mfc$ is $\mfm$-primary. 
If $A$ is $K_2$-regular, and $B/\mfm$ is separable over $A/\mfm$,
then $A = B.$ 
\end{prop}

\begin{proof}
Since $A$ is $K_2$-regular, it is $K_1$-regular. By Corollary
\ref{cor:K_1-regular_radical_conductor}, 
$\mfc = \mfm$, and $L = B/\mfm$ is 
a product of (separable, by hypothesis) field extensions of $K = A/\mfm$. 
It suffices to show that $K=L$.  
Let $A^\prime$ be the strict henselization of $A$ with respect to a
separable closure $K^\prime$ of $K$, $B^\prime = B\otimes_A A^\prime$,
and $L^\prime = L\otimes_K K^\prime = B^\prime/\mfm
B^\prime$. We will show that $K^\prime = L^\prime$ by
contradiction.

Since $B^\prime$ is a finite $A^\prime$-algebra, it is a product
$B^\prime = \prod_{i=1}^s B_i$ of local rings each with residue field
$K^\prime$. Write $\pi_i: B^\prime\to B_i$ for the projections.

Now assume for a contradiction that $s > 1$. Choose elements $f, \;
g\in\mfm$ such that $\pi_1(f)\neq 0$, $\pi_i(f) = 0$ for all
$i \neq 1$, $\pi_2(g)\neq 0$, and $\pi_i(g) = 0$ for all $i\neq 2$,
and $f$ and $g$ are linearly independent modulo
$\mfm^2$. Since the elements $f,\; gt\in A^\prime[t]$
satisfy $fgt = 0$, the Dennis-Stein symbol 
$\langle f,gt\rangle\in NK_2(A^\prime)$ is defined (see \cite{DS73} and
\cite[Theorem III.5.11.1]{WK}). 

The Dennis trace $K_2(A'[t])\to \Omega^2_{A^\prime [t]}$
sends $\langle f,gt\rangle\in NK_2(A^\prime)$ to
$-df\wedge d(gt)$ in $\Omega^2_{A^\prime [t]}$;
see e.g.\ \cite[2.2]{GW94}.
Because $f$ and $g$ are linearly
independent modulo $\mfm^2$, $df\wedge d(gt)\neq 0$ in
$\Omega^2_{A^\prime [t]}$, and it follows  
that the Dennis trace of
$\langle f,gt\rangle\in NK_2(A^\prime)$ is non-zero.  
Therefore $NK_2(A^\prime)\neq 0$.

On the other hand, it follows from \cite[Theorem 3.2]{vdK86} that
$NK_2(A) = 0$ implies $NK_2(A^\prime) = 0$. Thus, we have arrived at a
contradiction.	
\end{proof}

\begin{thm}\label{thm:K_2-regular-normal}
Let $A$ be a noetherian commutative ring containing a field 
and with finite normalization. 
Assume either that $\mathbb{Q}\subseteq A$, 
or that $\mathrm{dim}(A) = 2$ and $A/\mathfrak{m}$ is perfect 
for all maximal ideals $\mathfrak{m}$. If
$A$ is $K_2$-regular, then $A$ is normal.
\end{thm}

\begin{proof}
Assume first that $\mathbb{Q}\subseteq A$. 
We proceed by noetherian induction. Since $A$ is $K_2$-regular, it is
$K_1$-regular and hence reduced. It follows that the localization of
$A$ at any minimal prime is a field, and hence normal.
Now let $\mathfrak{p}$ be a prime
ideal such that $A_\mathfrak{q}$ is normal for all primes
$\mathfrak{q}$ properly contained in $\mathfrak{p}$. Then the local
ring $A_\mathfrak{p}$ is $K_2$-regular by Vorst \cite[1.9]{Vorst1},
and satisfies the hypotheses of Proposition \ref{prop:isolated-normal}.
We conclude that $A_\mathfrak{p}$ is normal. 
It follows by induction that all local rings of $A$ are normal,
and hence that $A$ is normal.

Now suppose that $\mathrm{dim}(A) = 2$. By \cite{Vorst1}, $A$ is
regular in codimension $1$. It follows that the conductor of the
normalization is $\mathfrak{m}$-primary. If, moreover, the residue
field is perfect, then we can again apply Proposition
\ref{prop:isolated-normal} to finish the proof.
\end{proof}

\newpage

\section{Higher du Bois singularities, $K$-theory and regularity.}
\label{sec:duBois}

 In this section we point out that the proof of Vorst's conjecture for
 an affine scheme $X$ of finite type over a field of characteristic
 zero contained in \cite{chw-vorst} does not actually require the full
 strength of the $K$-regularity hypothesis; instead, it suffices to
 assume only an algebraic-geometric consequence of this hypothesis,
 namely, that $X$ have {\it higher du Bois singularities}.

We begin by briefly recalling (in a form convenient for our
purposes) the definition of the du Bois complexes, $\underline{\Omega}^p_X$,
and the notion of du Bois and higher du Bois singularities. 
Du Bois complexes were introduced in \cite{duBois} to study 
mixed Hodge structures on possibly singular
quasiprojective complex varieties, and
du Bois singularities were introduced by Steenbrink in
\cite{Steenbrink83} as a weakening of the notion of rational singularities better behaved in families. 

Throughout this section, $k$ will be a field of characteristic zero,
and $\mathrm{Sch}/k$ will be the category of separated schemes
essentially of finite type over $k$. Following \cite{chww-bass2}, we
write $a:
\left(\mathrm{Sch}/k\right)_\cdh\to\left(\mathrm{Sch}/k\right)_\zar$ 
for the morphism of sites given by comparing the $cdh$ and Zariski topologies. 

\begin{defn}\label{defn:duBois-complex}
Let $X\in \mathrm{Sch}/k$, and $p\geq 0$. The $p$-th du Bois complex
of $X$ is the complex 
$\underline{\Omega}^p_{X/k} =\R a_* a^*\Omega^p_{/k} \vert_X$ 
of Zariski sheaves. 

Unless it is necessary to emphasize the base field, 
we will suppress it from the notation and simply write $\Omega^p_X$ 
and $\underline{\Omega}^p_{X}$ for differentials and 
du Bois complexes over $k$, respectively. 

A scheme $X$ is said to have {\it du Bois singularities} 
if $\cO_X\to \underline{\Omega}^0_X$ is an isomorphism.
\end{defn}

\begin{rem}\label{rem:duBois}
\begin{enumerate}
\item[(a)] This
definition of $\underline{\Omega}^p_X$
agrees with the original one given by du Bois in
  \cite{duBois} (using hyper-resolutions); 
see \cite[Lemma 2.1]{chww-bass2} for a proof of this fact.
\item[(b)] By adjunction, there is a natural homomorphism
  $\Omega^p_X\to \underline{\Omega}^p_X$. This homomorphism is an
  isomorphism when $X$ is smooth.
\item[(c)] The complex of sheaves of $\cO_X$-modules 
  $\underline{\Omega}^p_X$ has coherent cohomology sheaves 
(see \cite[Propn.\;4.4]{duBois}).
	\end{enumerate}
\end{rem}

The following generalization of the definition of du Bois singularity
has been studied recently, see \cite{JKSY22}, \cite{MP22}, 
\cite{SVV23}, \cite{FL24a}, \cite{FL24b}.

\begin{defn}\label{defn:dB-sing}
Let $X$ be a $k$-scheme, and $p\geq 0$. We say $X$ 
has (only) {\it $p$-\duBois singularities} if the natural homomorphisms
\[ \Omega^s_X\longrightarrow \underline{\Omega}^s_X\]
are isomorphisms for all $s\leq p$. 
\end{defn}

\begin{rem}
As pointed out  
in the Introduction of \cite{SVV23}, 
this notion (called ``strict $p$-du Bois" in loc.\,cit.) 
is poorly understood for $p>0$, 
unless $X$ is a local complete intersection. 
We will return to this in section \ref{sec:lci}.
\end{rem}

We will connect du Bois singularities to $K$-regularity using descent for Hochschild 
homology. Recall that  $HH(/k)$ satisfies $p$-cdh descent on 
a $k$-scheme $X$ if the fiber of the map from $HH(X/k)$ to its
$cdh$-hypercohomology is $p$-connected.

\begin{lem}\cite[Lemma 2.3]{chw-vorst}\label{lem:tech1}
Let $R$ be a $k$-algebra essentially of finite type, and let $p\geq 0$. Then $HH(/k)$ satisfies $p$-cdh
descent on $X=\Spec(R)$ if and only if the following three
conditions hold simultaneously:
\addtocounter{equation}{-1}
\begin{subequations}
\begin{gather}
HH_m^{(q)}(R/k)=0 \text{\quad if } 0\le q<m\le p;\label{hhkvanish}
\\
\Omega^q_{R/k}\rightarrow H^0_{\cdh}(X,\Omega^q_{/k})\;
\text{is bijective if } q\le p, \text{ onto if } q=p+1;\label{agreek}
\\
H^s_{\cdh}(X,\Omega^q_{/k})=0 \quad\text{ if } s>0 \text{ and }0\le q\le s+p+1.
\label{cdhkvanish}
\end{gather}
\end{subequations}
\end{lem}

The authors and Corti\~nas proved in \cite[Theorem 3.1]{chw-vorst}
that a finite type $k$-algebra satisfying the equivalent conditions of
Lemma \ref{lem:tech1} is regular in codimension $p$. The following
theorem shows that the hypothesis can be weakened:

\begin{thm}\label{thm:dB-Vorst}
Let $X\in\mathrm{Sch}/k$. If $X$ has $p$-\duBois singularities, then
$X$ is regular in codimension $p$.
\end{thm}

\begin{proof}
Let $x\in X$ be a point of codimension $d \leq p$. 
By \cite[Lemma 4.10]{chw-vorst}, we can represent the local ring $A = \cO_{X,x}$
as the localization of a finite type $F$-algebra 
at a maximal ideal, for some field $F$ containing $k$. 
The proofs of \cite[Lemma 4.4 and Prop.\;4.8]{chw-vorst} go through,
noting that the hypotheses of Theorem \ref{thm:dB-Vorst} suffice, showing that 
if $X$ is $p$-\duBois (over $k$) then $\Spec(A)$ is 
$p$-\duBois (over $F$). Thus we may assume,
without loss of generality, that $X = \Spec (A)$ is the spectrum of a local ring of 
a $d$-dimensional $k$-algebra of finite type at a maximal ideal, and $x$ is its closed point.
The hypothesis that $X$ has $p$-\duBois singularities then 
is equivalent to the following subset
 of the conditions of Lemma \ref{lem:tech1} (recall that $d\leq p$): 
\addtocounter{equation}{-1}
\begin{subequations}
	\begin{gather}
		\Omega^q_{A/k}\rightarrow H^0_{\cdh}(X,\Omega^q_{/k})\quad
		\text{is an isomorphism if } q\le d;
		\\
		H^s_{\cdh}(X,\Omega^q_{/k})=0 \quad\text{ if } s>0 \text{ and }0\le q\le d.
	\end{gather}
\end{subequations}
By inspection, the proof of \cite[Theorem 3.1]{chw-vorst} only needs
the weaker hypotheses given here, not all of Lemma 2.3 loc.\;cit. 
Indeed, the properties above suffice to conclude, 
as in the proof of \cite[Theorem 3.1]{chw-vorst}, that
\[
H^*(\Omega^{\le d}_{A/k},d)\cong 
\bbHcdh^{*}(X,\Omega^{\le d}_{/k}) \cong 
\bbHcdh^{*}(X,\Omega^{\bullet}_{/k})
\]
and, as explained there, this implies that $\Omega^{d+1}_{A/k} = 0$. 
 This proves the assertion.
\end{proof}

The following is extracted from \cite{chw-vorst}.

\begin{prop}\label{prop:K-regular-to-duBois}
Let $X\in \mathrm{Sch}/k$ be an affine scheme. If $X$ is
$K_{p+1}$-regular, then $X$ has $p$-\duBois singularities.
\end{prop}

\begin{proof}
By \cite[Corollary 1.7]{chw-vorst}, the $K$-regularity hypothesis is
equivalent to the assertion that $HH(/\mathbb{Q})$ satisfies
$p$-cdh-descent on $X$.
Now \cite[Propn.\,4.8]{chw-vorst} implies that $HH(/k)$ satisfies
$p$-cdh-descent on $X$ as well. Finally Lemma \ref{lem:tech1} gives an
explicit description of this descent property; parts (\ref{agreek})
and (\ref{cdhkvanish}) of this description in particular imply that
$X$ has $p$-\duBois singularities.
\end{proof}

Combining the above yields 
a new proof of the main theorem of \cite{chw-vorst}.
Recall that $k$ is a field of characteristic~0.

\begin{thm}\label{thm:Vorst-redux}
Let $X\in\mathrm{Sch}/k$ be an affine scheme, and $p\geq 1$. 
If $X$ is $K_{p+1}$-regular, then $X$ is 
regular in codimension $p$.
\end{thm}

\begin{proof}
Proposition \ref{prop:K-regular-to-duBois} asserts that $X$ is $p$-\duBois. It now follows from Theorem \ref{thm:dB-Vorst} that $X$ is regular in codimension $p$. 
\end{proof}

\section{The case of local complete intersections.}\label{sec:lci}

In this short section, we combine Hodge theoretic results of 
Musta\c t\u a and Popa \cite{MP22} 
with the previous section to prove that the $K$-regularity bound
for singularities given by Vorst's conjecture and 
proven in \cite{chw-vorst} can be significantly
improved if the singularity is a local complete intersection ($lci$).
We also formulate some questions.

The following theorem is proved by Musta\c t\u a and Popa in \cite{MP22} 
and is stated in the form used here (with $k = \mathbb{C}$) 
as \cite[Theorem 3.8]{FL24b}. The hypersurface case of this theorem also follows from bounds proved in \cite{MOPW23} and \cite{JKSY22}.

\begin{thm}[Musta\c t\u a--Popa]\label{thm:dB-smooth}
Let $X/k$ be locally a complete intersection 
with $p$-\duBois singularities, $p\geq 1$. Then 
$X$ is normal and regular in codimension $2p$.  
\end{thm}

\begin{question}
Let $A$ be a $k$-algebra, essentially of finite type over $k$.
As can be seen using Lemma \ref{lem:tech1} 
and \cite[Propn.\;4.8]{chw-vorst},   $A$ is $K_2$-regular 
if, and only if, $\Spec (A)$ is $1$-\duBois, the natural homomorphism 
$\Omega^2_{A/k}\rightarrow H^0_{\cdh}(X,\Omega^2_{/k})$ is onto, and 
certain cohomology groups of the du Bois complexes vanish. 
We proved in Theorem \ref{thm:K_2-regular-normal} that 
a $K_2$-regular algebra is normal, without the assumption 
that $A$ is a local complete intersection. 

This suggests the following question: Suppose $A$ is a $1$-\duBois
algebra, not necessarily a local complete intersection. 
Is $A$ normal? We note that $A$ is
regular in codimension $1$ by Theorem \ref{thm:Vorst-redux}, so this
question really asks if $A$ necessarily satisfies Serre's condition (S2).
\end{question}

The following result provides the promised improvement of Vorst's bound
in the case of a local complete intersection.

\begin{thm}\label{thm:Vorst-lci}
Let $X\in\mathrm{Sch}/k$ be an affine scheme, and $p\geq 1$. 
If $X$ is $K_{p+1}$-regular and a local complete intersection, 
then $X$ is normal and regular in codimension $2p$.
\end{thm}

\begin{proof}
This follows immediately from Proposition \ref{prop:K-regular-to-duBois} 
and Theorem \ref{thm:dB-smooth}.
\end{proof}

\begin{question} 
We wonder if the
lci hypothesis is needed to prove the sharper bound. 
Typically, schemes that are not local complete intersection are 
not even $1$-\duBois,
since their sheaves of differentials are usually not torsion free. 
On the other hand, the arguments developed in \cite{MP22} employ the
Hodge structure on (top) local cohomology and rely crucially on the
lci condition. 

The first case of interest is those of surfaces, which suggests the
question: Suppose $X$ is an affine surface over a field of
characteristic zero. If $X$ is $K_2$-regular, is it necessarily
regular?
\end{question}

Theorem \ref{thm:dB-smooth} can be sharpened; the following was communicated to us by Wanchun Shen. Let $S$ be a smooth $k$-scheme, and $X\subseteq S$ a local complete intersection subvariety. The minimal exponent $\tilde{\alpha}(S,X)$ is a positive real number; it was introduced for hypersurfaces in \cite{Saito94}, and defined for general local complete intersection subschemes in \cite{CDMO24}. 

Given a complete intersection local ring $A$, the \emph{minimal embedding codimension} of $A$ is the minimal integer $r$ such that there exists a regular local ring $R$ and regular sequence $\left(f_1, \dotsc, f_r \right)$ in $R$ such that $A\cong R/\left(f_1, \dotsc, f_r \right)$. 

Combining results from \cite{CDMO24} and \cite{CDM24} we obtain the following result:

\begin{thm}[Wanchun Shen]\label{thm:dB-codim}
Let $A = \cO_{X,x}$ be a local ring of a local complete intersection $k$-variety at a closed point, and $r$ the minimal embedding codimension of $A.$
Suppose $X  = \Spec (A)$ has $p$-\duBois singularities, $p\geq 1$. Then $X$ is regular in codimension $2p + r - 1$.
\end{thm}

\begin{proof}
We can find a smooth $k$-variety $S$ and closed point $s\in S$ such that $A\cong \cO_{S,s}/J$ for an ideal $J$ generated by a regular sequence of length $r$. Let $\mathfrak{m}$ be the maximal ideal of $A$. Because $r$ is minimal, $J\subseteq \mathfrak{m}^2$. 
Intersecting with generic hyperplanes preserves the $p$-\duBois property by \cite[Corollary 5.11]{FL24b}, so we may assume that $X$ has an isolated singularity. Let $n = \dim S$ and $d = n - r = \dim X$. Since $J \subseteq \mathfrak{m}^2$, \cite[Theorem 1.2 (iii)]{CDMO24} implies that $\tilde{\alpha}(S,X)\leq n/2 = (d+r)/2.$ On the other hand, because $X$ is $p$-du\;Bois, \cite[Theorem F]{MP22} shows that $\tilde{\alpha}(S,X)\geq p+r.$ (The minimal exponent is not used in {\em loc.\,cit.}, but it is explained in \cite[Subsection 2.4]{CDM24} that Theorem F can be interpreted as cited.)

Combining the two inequalities, we get that
\[2(p+r) \leq d+r \Longrightarrow 2p + r \leq d.\]

Since $X$ has an isolated singularity, we conclude that it is regular in codimension $2p + r - 1$, as asserted.
\end{proof}

\begin{cor}\label{cor:dB-codim}
	Let $(X,x) = \left(\Spec (A),\mathfrak{m}\right)$ be a local, complete intersection singularity of dimension $d$ over $k$, and let $r$ be its minimal embedding codimension. If $X$ has $p$-du Bois singularities, $p\geq 1$, then $r \leq d - 2p$.	
	\end{cor}
	
	\begin{proof}
	 Assume for a contradiction that the minimal embedding codimension is $r > d-2p$. Then $X$ is regular in codimension $d$, and hence, regular. But this implies $r = 0$, a contradiction.  	
	\end{proof}

\begin{rem}\label{rem:Vorst-lci}
	We can apply Corollary \ref{cor:dB-codim} to obtain $K$-regularity bounds on the embedding codimension of germs of singularities. For example, if $(X,x)$ is a local, complete intersection singularity of
	 dimension $3$, and $X$ is $K_2$-regular, then $X$ is a hypersurface. 
\end{rem}

\medskip
\subsection*{Acknowledgements}
The authors thank Wanchun Shen for generously sharing Theorem \ref{thm:dB-codim} with them. Haesemeyer would like to thank the Australian Research Council for supporting this work through DP--210103397. Weibel would like to thank the National Science Foundation for support on this project through the grant 2001417, and the University of Melbourne for hosting him during work on this paper. 

\smallskip
\subsection*{Data availability} No data were generated or used. 

\smallskip
\subsection*{Conflict of interest declaration} On behalf of all authors, the corresponding author states that there is no conflict of interest.


\end{document}